\documentclass{biometrika}

\usepackage{amsmath}

\usepackage{times}
\usepackage{bm}
\usepackage{natbib}

\usepackage[plain,noend]{algorithm2e}
\usepackage{latexsym,amsmath,amssymb,amsfonts}

\makeatletter
\renewcommand{\algocf@captiontext}[2]{#1\algocf@typo. \AlCapFnt{}#2} 
\def\@algocf@capt@plain{top}
\renewcommand{\algocf@makecaption}[2]{%
  \addtolength{\hsize}{\algomargin}%
  \sbox\@tempboxa{\algocf@captiontext{#1}{#2}}%
  \ifdim\wd\@tempboxa >\hsize
    \hskip .5\algomargin%
    \parbox[t]{\hsize}{\algocf@captiontext{#1}{#2}}
  \else%
    \global\@minipagefalse%
    \hbox to\hsize{\box\@tempboxa}
  \fi%
  \addtolength{\hsize}{-\algomargin}%
}
\makeatother

\begin{document}

\jname{Accepted for publication in Biometrika published by Oxford University Press}


\markboth{T.~B.~Berrett}{Discussion}

\title{Discussion of `Multiscale Fisher's independence test for multivariate dependence'}

\author{T.~B.~BERRETT}
\affil{Department of Statistics,  University of Warwick,\\ Coventry CV4 7AL, U.K.
\email{tom.berrett@warwick.ac.uk}}

\maketitle

\section{Introduction}

We would like to congratulate the authors on an interesting paper that provides many avenues for discussion. The reduction of a nonparametric problem to a sequence of $2 \times 2$ contingency tables is very elegant, and Fisher's exact test is used to great effect. There are now a large number of methods for multivariate independence testing, but the proposed procedure stands out by providing uniform Type I error controls without a large computational burden, combining the speed of an asymptotic test with the rigorous validity of a permutation test. Moreover, the ability of the test to provide information on the form of the possible dependence is a feature that should be appreciated in practice, and one that is worthy of further consideration.

The problem of independence testing has received a great deal of attention throughout the history of statistics, and is a canonical example of a nonparametric testing problem where the null hypothesis is composite. Further, it is a problem where departures from the null hypothesis can take many forms. Beyond the obvious practical importance of independence testing, these two features of the problem are key reasons for the development of such a vast literature on the topic. In this discussion we describe how the new work fits in to this literature, and outline criteria to consider when choosing an independence test.

\section{Validity of tests and computation}

When null hypotheses are simple the null distribution of any test statistic is, in principle, known, and critical values can be tabulated. However, the composite nature of the null hypothesis in independence testing means that ensuring the validity of a test is non-trivial. There are at least three main approaches taken in the literature. First, with any test statistic it is possible to resample the data to carry out a permutation test, to provide exact Type I error control with a heavier computational cost \citep[e.g.][]{berrett2019nonparametric,berrett2020optimal,berrett2021usp}. One way to reduce this cost is to use a test statistic whose distribution under the null hypothesis can be asymptotically approximated by a known distribution, and there are consequently many works on independence testing that derive an asymptotic null distribution as a main result \citep[e.g.][]{szekely2007measuring,gretton2007kernel}. 

The third approach, taken in the current work, is to carefully design test statistics that are distribution-free, in that their distribution is the same for all elements of the null hypothesis space. The construction of such test statistics when marginal distributions are univariate is relatively straightforward, as one can apply the empirical distribution function to the data and work with ranks, and classical methods such as Kendall's $\tau$, Spearman's rank correlation coefficient, and Hoeffding's $D$ take advantage of this fact. With multivariate data this is much more difficult, though there has been notable success in using optimal transport to design appropriate multivariate ranks \citep{chernozhukov2017monge}, and this method has been used successfully in independence testing \citep{deb2021multivariate,shi2021distribution}. By applying a suitable quantile transformation to the data the marginal distributions are standardized, effectively reducing the composite null hypothesis to a point. The authors of the current work remark that their test is conditional on the marginal values, a fact which is inherited from Fisher's exact test.

Where the current proposal differs from previous uniformly valid multivariate approaches is in computational complexity. With permutation tests any test statistic must be calculated $B$ times, for some $B \in \mathbb{N}$ that is typically in the hundreds. The time taken to compute optimal transport-based multivariate ranks is $O(n^3)$ in the worst case, where $n$ is the sample size \citep{deb2021multivariate}. By contrast, the proposed test at a fixed resolution can be carried out in $O(n \log n)$ time, potentially providing a marked improvement. It is interesting to contrast this computational cost with that of permutation $k$-nearest neighbour methods  \citep[e.g.][]{berrett2019nonparametric}, which run in $O(Bkn \log n)$ time. If we are interested in having high power against alternatives with local forms of dependence then the current method should be used with a fine resolution, increasing the computational cost, while a $k$-nearest neighbour method should be used with a small value of $k$, decreasing the computational cost. As the worst-case computational time of the new method appears to be exponentially increasing in the resolution, it may be the case that such local forms of dependence can be detected more quickly with alternative methods.

\section{Power results}

Controlling the Type I error rate of a test is, of course, only part of the problem, and designing tests with high power is more complicated still. In Theorem~3 the proposed test is shown to be consistent, in that its power converges to one as the sample size increases and the resolution becomes finer, for each fixed distribution that does not satisfy the null hypothesis. This is a very desirable property for an independence test, but the power does not converge to one uniformly over alternative distributions. While it is possible, for example with permutation tests, to control the Type I error uniformly over the null hypothesis space, when data is continuous it is impossible to design an independence test that has power against all forms of dependence simultaneously, even after fixing the strength of the dependence \citep{zhang2019bet,berrett2020optimal}. This is another reason for the great many different test statistics in the literature, each of them making implicit or explicit choices about which kinds of dependence are prioritized for detection. The authors assume that the maximal resolution at which exhaustive testing is applied $R^* = o(\log n)$, so the current method naturally focusses on dependence occurring at the coarsest resolutions, progressing to the finer scales only if enough data is available.

As an illustration of this phenomenon, consider distributions exhibiting sinusoidal dependence \citep{sejdinovic2013equivalence,berrett2019nonparametric} with density functions
\[
	f_{\rho,L}(x,y) = 1 + \rho \sin(2^{L+1} \pi x) \sin(2^{L+1} \pi x)
\]
on $[0,1]^2$, for $\rho \in [0,1)$ and $L \in \mathbb{N}$. Here $\rho$ controls the strength of the dependence, and it can be seen that the mutual information between $X$ and $Y$ when $(X,Y)$ has density $f_{\rho,L}$ is a function of $\rho$ only. On the other hand, $L$ controls the frequency at which dependence occurs, controlling the form of the dependence. For any $(k_1,k_2) \in \mathbb{N}^2 \setminus \{L+1,L+2,\ldots\}$, any $\ell_1 \in \{1,\ldots,2^{k_1}\}$ and any $\ell_2 \in \{1,\ldots,2^{k_2}\}$, we have
\begin{align*}
	&\int_{(\ell_1-1)2^{-k_1}}^{\ell_1 2^{-k_1}} \int_{(\ell_2-1)2^{-k_2}}^{\ell_2 2^{-k_2}} f_{\rho,L}(x,y) \,dy \,dx \\
	&= 2^{-k_1-k_2} \int_0^1 \int_0^1 \bigl\{ 1 + \rho \sin(2^{L-k_1+1} \pi u) \sin(2^{L-k_2+1} \pi v) \bigr\} \,dv \,du = 2^{-k_1-k_2},
\end{align*}
and we can therefore see that the proposed test can only detect dependence if the chosen resolution exceeds $2L$. Naturally, the test has best power against dependence at low frequencies, as do most tests used in practice, but this is a methodological choice and one could design a test targeting any frequencies that were of most interest \citep{berrett2020optimal}.

In order to compare the power of different valid independence tests we must create alternative hypotheses, stratified by both strength and form of dependence, that reflect the types of alternatives we are most interested in. Studying power is thus necessarily more subjective than studying validity, but once we have formulated our assumptions and priorities it is possible to search for tests with optimal power. The most common approaches to nonparametric statistics impose smoothness conditions, for example Sobolev smoothness, and prioritize alternative distributions that are both smooth and exhibit relatively strong dependence \citep{berrett2020optimal,albert2022adaptive}. However, there are other choices: \cite{zhang2019bet} imposes conditions to assume that alternatives exhibit dependence early on in the binary expansions of the data, and the current work assumes, for example in Theorem~3, that alternatives exhibit dependence at relatively coarse resolutions. These three types of assumptions are somewhat related, especially for the alternatives $f_{\rho,L}$ above, but there is no guarantee that optimal tests with respect to one set of assumptions are optimal with respect to another. 

\section{Providing intuition about form of dependence}

Independence tests do not typically give any information to the practitioner beyond the value of the test statistic and the decision to reject the null hypothesis or not. Due to the large number of possible departures from the null hypothesis, it is very difficult to provide useful information on the form of the dependence in the data. The sequential nature of the new test naturally allows for more information to be given with a rejection: as the procedure searches for subsets of the sample space that exhibit dependence, it is informative to report where this dependence was detected. We remark that \citet{zhang2019bet} defines another procedure giving interpretable results, which works by finding the coefficients in the binary expansions of the data that are maximally correlated. Despite this two approaches, this aspect of the independence testing problem is rarely discussed in the literature.

It would be interesting to see more research into this topic in the future. As with the design of powerful tests, it seems that there should be many approaches whose suitability depends on the forms of dependence of most interest. The current proposal will typically return a small subset of the sample space where there is evidence of dependence. For the sinusoidal densities $f_{\rho,L}$ described above, it seems likely to return a set of the form $[(\ell_1-1)2^{-k_1},\ell_12^{-k_1}) \times [(\ell_2-1)2^{-k_2},\ell_22^{-k_2})$, where $k_1,k_2 \geq L+1$. This very local view does not tell the full story for this particular setting, and an approach such as \cite{zhang2019bet} may be more successful.

More generally, the problem of summarizing a potentially complex relationship between two variables is very interesting. Writing $f$ for the density of $(X,Y)$ and $f_X$ and $f_Y$ for the marginal densities of $X$ and $Y$, respectively, we may aim to return an approximation of the set $\{(x,y) : f(x,y) > f_X(x)f_Y(y) \}$ where the joint density is larger than the product of the marginal densities. Alternatively, motivated by the maximal correlation advocated by \cite{renyi1959measures}, we may aim to find functions $F$ and $G$ within some interpretable class to maximize the correlation between $F(X)$ and $G(Y)$. Both suggestions would present significant statistical challenges, and would likely involve modelling choices, but could be worth exploration.

\section{Conclusion}

Despite its long history, independence testing remains an active and exciting area of research where many delicate statistical issues can be debated. In this discussion we have tried to provide an overview of the topic of independence testing, particularly focussing on the parts that are relevant to the current submission. The primary concern in the problem is to ensure that tests control the rate of Type I error, and can therefore provide the basis for valid inference. The authors provide a neat way of doing this, based on elegant mathematical arguments, that bypasses the computational cost associated with many valid independence tests. Maximizing the power of a test is a more subtle issue, but the authors show that their method is consistent against all fixed alternatives. One great advantage of the method, in our opinion, is its interpretability, and we feel that this deserves more attention from the community.


\begin{thebibliography}{7}
\expandafter\ifx\csname natexlab\endcsname\relax\def\natexlab#1{#1}\fi

\bibitem[{Albert et al.(2022)}]{albert2022adaptive}
\textsc{Albert, M., Laurent, B., Marrel, A. \& Meynaoui, A.} (2022).
\newblock Adaptive test of independence based on HSIC measures
\newblock \textit{Ann. Statist., to appear}.

\bibitem[{Berrett and Samworth(2019)}]{berrett2019nonparametric}
\textsc{Berrett, T. B. \& Samworth, R. J.} (2019).
\newblock Nonparametric independence testing via mutual information
\newblock \textit{Biometrika} \textbf{106}, 547--566.

\bibitem[{Berrett et al.(2021)}]{berrett2020optimal}
\textsc{Berrett, T. B., Kontoyiannis, I. \& Samworth, R. J.} (2021).
\newblock Optimal rates for independence testing via $U$-statistic permutation tests
\newblock \textit{Ann. Statist.} \textbf{49}, 2457--2490.

\bibitem[{Berrett and Samworth(2021)}]{berrett2021usp}
\textsc{Berrett, T. B. \& Samworth, R. J.} (2021).
\newblock USP: an independence test that improves on Pearson's chi-squared and the $G$-test
\newblock \textit{Proc. R. Soc. A.} \textbf{477}: 20210549.

\bibitem[{Chernozhukov et al.(2017)}]{chernozhukov2017monge}
\textsc{Chernozhukov, V., Galichon, A., Hallin, M. \& Henry, M.} (2017).
\newblock Monge--Kantorovich depth, quantiles, ranks and signs
\newblock \textit{Ann. Statist.} \textbf{45}, 223--256.

\bibitem[{Deb and Sen(2022)}]{deb2021multivariate}
\textsc{Deb, N. \& Sen, B.} (2022).
\newblock Multivariate rank-based distribution-free nonparametric testing using measure transportation
\newblock \textit{J. Amer. Statist. Assoc., to appear.}

\bibitem[{Gretton et al.(2007)}]{gretton2007kernel}
\textsc{Gretton, A., Fukumizu, K., Teo, C., Song, L., Sch\"olkopf, B. \& Smola, A.} (2007).
\newblock A kernel statistical test of independence
\newblock \textit{NeurIPS 20}

\bibitem[{R\'enyi(1959)}]{renyi1959measures}
\textsc{R\'enyi, A.} (1959).
\newblock On measures of dependence
\newblock \textit{Acta Math. Hungar.} \textbf{10}, 441--451.

\bibitem[{Sejdinovic et al.(2013)}]{sejdinovic2013equivalence}
\textsc{Sejdinovic, D., Sriperumbudur, B., Gretton, A. \& Fukumizu, K.} (2013).
\newblock Equivalence of distance-based and RKHS-based statistics in hypothesis testing
\newblock \textit{Ann. Statist.} \textbf{41}, 2263--2291.

\bibitem[{Shi et al.(2022)}]{shi2021distribution}
\textsc{Shi, H., Drton, M. \& Han, F.} (2022).
\newblock Distribution-free consistent independence tests via center-outward ranks and signs
\newblock \textit{J. Amer. Statist. Assoc.} \textbf{117}, 395--410.

\bibitem[{Sz\'ekely et al.(2007)}]{szekely2007measuring}
\textsc{Sz\'ekely, G. J., Rizzo, M. L. \& Bakirov, N. K.} (2007).
\newblock Measuring and testing dependence by correlation of distances
\newblock \textit{Ann. Statist.} \textbf{35}, 2769--2794.

\bibitem[{Zhang(2019)}]{zhang2019bet}
\textsc{Zhang, K.} (2019).
\newblock BET on independence
\newblock \textit{J. Amer. Statist. Assoc.} \textbf{114}, 1620--1637.

\end{thebibliography}
\end{document}